\documentclass[reqno, 12]{amsart}
\usepackage{amsmath}
\usepackage{amscd,amsthm,amsfonts,amsopn,amssymb,verbatim}
\parskip =\medskipamount
\numberwithin{equation}{section}
\newtheorem{theorem}{Theorem} 

\newtheorem{proposition}[theorem]{Proposition}

\theoremstyle{remark}
\DeclareMathOperator{\supp}{supp\,}

\def\be{\begin{equation}}
\def\ee{\end{equation}}
\def\vp{\varphi}

\def\mod{\text{\rm mod\,}}
\allowdisplaybreaks

\begin{document}

\title{A note on the Schr\"odinger maximal function}
\date{\today}
\author{J.~Bourgain}
\address{J.~Bourgain, Institute for Advanced Study, Princeton, NJ 08540}
\email{bourgain@math.ias.edu}
\thanks{The author was partially supported by NSF grants DMS-1301619}
\begin{abstract}
It is shown that control of the Schr\"odinger maximal function $\sup_{0<t<1} |e^{it\Delta} f|$ for $f\in H^s(\mathbb R^n)$ requires $s\geq \frac n{2(n+1)}$
\end{abstract}

\maketitle

\section
{Introduction}

Recall that the solution of the linear Schr\"odinger equation
\be\label{1.1}
\begin{cases}
iu_t-\Delta u=0\\
u(x, 0)= f(x)
\end{cases}
\ee
with $(x, t)\in\mathbb R^n\times\mathbb R$ is given by
\be\label {1.2}
e^{it\Delta} f(x) =(2\pi)^{-n/2}\int e^{i(x.\xi+t|\xi|^2)}\hat f(\xi)d\xi.
\ee
Assuming $f$ belongs to the space $H^s(\mathbb R^n)$ for suitable $s$, when does the almost convergence property
\be\label{1.3}
\lim_{t\to 0} e^{it\Delta} f=f \ \ a.e.
\ee
hold?
The problem was brought up in Carleson's paper \cite{C} who proved convergence for $s\geq \frac 14 $ when $n=1$.
Dahlberg and Kenig \cite {D-K} showed that this result is sharp.
In higher dimension, the question of identifying the optimal exponent $s$ has been studied by several authors and our state of
knowledge may be summarized as follows.
For $n=2$, the strongest result to date appears in \cite{L} and asserts \eqref{1.3} for $f\in H^s (\mathbb R^2)$, $s>\frac 38$.
More generally, for $n\geq 2$ \eqref{1.3} was shown to hold for $f\in H^s(\mathbb R^n)$, $s>\frac {2n-1}{4n}$ (see \cite {B}).

In the opposite direction, for $n\geq 2$ the condition $s\geq \frac n{2(n+2)}$ was proven to be necessary (see \cite{L-R} and also \cite{D-G} for a different approach
based on pseudoconformal transformation). 
Here we show the following stronger statement.

\begin{proposition}\label{Proposition1}

Let $n\geq 2$ and $s< \frac n{2(n+1)}$.
Then there exist $R_k\to\infty$ and $f_k\in L^2(\mathbb R^n)$ with $\hat f_k$ supported in the annulus $|\xi|\sim R_k$, such that $\Vert f_k\Vert_2=1$ and
\be\label{1.4}
\lim_{k\to\infty} R_k^{-s}\Big\Vert\sup_{0<t<1}| e^{it\Delta} f_k (x)|\big\Vert_{L^1(B(0, 1))} =\infty.
\ee
\end{proposition}

There is some evidence the exponent $\frac n{2(n+1)}$ could be the optimal one, though limited to multi-linear considerations appearing in \cite {B}.
Of course, the $n=1$ case coincides  with the \cite {D-K} result, while for $n=2$, the above Proposition leaves a gap between $\frac 13$ and $\frac 38$.
It may be also worth to point out that for $n=2$, in some sense, our example fits a scenario where the arguments from \cite {B} require the $s>\frac 38$ condition.

\section
{An example}

Denote $x=(x_1, \ldots, x_n)=(x_1, x')\in B(0, 1)\subset \mathbb R^n$.
Let $\vp:\mathbb R\to\mathbb R_+, \Phi: \mathbb R^{n-1}\to \mathbb R_+$ satisfy $\supp \hat\vp \subset [-1, 1]$, $\supp \hat\Phi \subset B(0, 1)$, $\hat \vp, \hat\Phi$ smooth and
$\vp(0)=\Phi(0)=1$.
Set $D=R^{\frac {n+2}{2(n+1)}}$ and define
\be\label{2.1}
f(x) =e(Rx_1)\vp(R^{\frac 12} x_1)\Phi(x') \prod^n_{j=2} \Big(\sum_{\frac R{2D}< \ell_j< \frac RD} e^{iD\ell_jx_j}\Big)
\ee
where $\ell =(\ell_2, \ldots, \ell_n)\in\mathbb Z^{n-1}$.
Hence
\be\label{2.2}
\Vert f\Vert_2 \sim R^{-\frac 14} \Big(\frac RD\Big)^{\frac {n-1} 2} \ \text { and } \ \supp\hat f\subset [|\xi|\sim R].
\ee
Clearly, denoting $e(z)=e^{iz}$,
$$
e^{it\Delta} f (x) =\iint \hat\vp (\lambda) \hat\Phi (\xi')\Big\{\sum_\ell e\big((R+\lambda R^{\frac 12})x_1+(\xi'+D\ell). x'+(R+\lambda R^{\frac 12})^2 t+|\xi'+D\ell|^2
t\big)\Big\} d\lambda d\xi'.
$$
Taking $|t|<\frac cR$, $|x|<c$ for suitable constant $c>0$, one gets
\begin{align}\label{2.3}
|e^{it\Delta} f(x)|&\sim \Big|\int \hat\vp(\lambda)\Big\{\sum_\ell e(\lambda R^{\frac 12} x_1 + D\ell.x'+2\lambda R^{\frac 32}t+D^2|\ell|^2t)\Big\}d\lambda\Big|\nonumber\\
&\sim \vp \big(R^{\frac 12} (x_1+2Rt)\big)\Big|\sum_\ell e(D\ell.x'+D^2|\ell|^2 t)\Big|
\end{align}

Specify further $t=- \frac {x_1}{2R}+\tau$ with $|\tau|<\frac  1{10} R^{-\frac 32}$ in order to ensure that the first factor in \eqref{2.3} should be $\sim 1$.
For this choice of $t$, the second factor becomes
\begin{align}\label{2.4}
&\Big|\sum_\ell e\Big( D\ell.x' -\frac {D^2}{2R} |\ell|^2 x_1+ D^2|\ell|^2 \tau\Big)\Big|=\nonumber\\
&\prod^n_{j=2} \Big|\sum_{\frac R{2D} <\ell_j<\frac RD} e\big(\ell_j y_j+\ell_j^2(y_1+s)\big)\Big|
\end{align}
with
\be\label{2.5}
y'=Dx'(\mod 2\pi) \qquad y_1=-\frac {D^2}{2R} x_1 (\mod 2\pi)
\ee
and where $s=D^2\tau$ is subject to the condition
\be\label{2.6}
|s|\lesssim D^2 R^{-3/2}=R^{-\frac{n-1}{2(n+1)}}.
\ee
We view $y=(y_1, y')$ as a point in the $n$-torus $\mathbb T^n$.
Next, define the following subset $\Omega\subset \mathbb T^n$
\be\label{2.7}
\Omega =\bigcup_{q\sim R^{\frac {n-1}{2(n+1)}}, a} \Big\{ (y_1, y'); \Big|y_1-2\pi \frac {a_1}q\Big|< cR^{-\frac {n-1}{2(n+1)}}
\ \text { and } \ \Big |y'-2\pi\frac {a'} q\Big|< c\frac DR\Big\}
\ee
with $a=(a_1, a')$ $(\mod q)$ and $(a_1, q)=1$.

Hence $|\Omega|\sim R^{\frac {n-1}{2(n+1)}} R^{n\frac {n-1}{2(n+1)}} R^{-\frac {n-1}{2(n+1)}}
\Big(\frac DR\Big)^{n-1} \sim 1$ and we take $x\in B(0, 1)$ for which $y$ given by \eqref {2.5} belongs to $\Omega$.
Clearly this gives a set of measure at least $c_1>0$.
We evaluate \eqref{2.4} for $y\in\Omega$.
Let $q\sim R^{\frac {n-1}{2(n+1)}}$ and $(a_1, a')$ $(\mod q)$ satisfy the approximations stated in \eqref{2.7} and set 
$s=2\pi\frac {a_1}q -y_1$ for which \eqref{2.6} holds.
Clearly for $j=2, \ldots, n$, by the quadratic Gauss sum evaluation
$$
\begin{aligned}
\Big|\sum_{\frac R{2D} <\ell_j< \frac RD} e\big(\ell_j y_j+\ell^2_j(y_1+s)\big)\Big|&\sim \Big|\sum _{\frac R{2D}<\ell_j<\frac RD}
e(2\pi \frac {a_j}q \ell_j+2\pi\frac {a_1}q \ell_j^2)\Big|\\
&\sim R^{\frac 1{2(n+1)}} \Big|\sum^{q-1}_{\ell_j=0}
e(2\pi \frac {a_j}q\ell_j+2\pi \frac {a_1}q \ell^2_j)\Big|\\
&\sim R^{\frac 1{2(n+1)}} q^{\frac 12}\sim R^{\frac 14}
\end{aligned}
$$
and
\be\label{2.8}
\eqref{2.4} \sim R^{\frac {n-1} 4}.
\ee
Recalling \eqref{2.2}, we obtain for $x\in B(0, 1)$ in a set of measure $c_1>0$
that
\be\label{2.9}
\sup_{0<t<1} \frac{|e^{it\Delta}f(x)|}{\Vert f\Vert_2} \gtrsim  R^{\frac {n-1}4} R^{\frac 14}\Big(\frac DR\Big)^{\frac {n-1}2}= R^{\frac n{2(n+1)}}
.\ee
The claim in the Proposition follows.

\end{document}